\documentclass[10pt]{article} %-*-LaTeX-*-
\usepackage[dvips]{graphicx}
\usepackage{psfrag}
\usepackage{calc}  %%required to make lists work correctly
\usepackage{amsmath}  
\usepackage{amsthm}
\usepackage{amsfonts}
\usepackage{amsopn}
\usepackage{amssymb}
\usepackage{amstext}
\usepackage[matrix,tips,graph,curve,web]{xy}

\makeatletter

\def\@maketitle{%
  \newpage
  \null
  \let \footnote \thanks
    {\normalfont\sffamily\bfseries\Large\noindent\@title \par}%
    \vskip 1em%
    {\normalfont\sffamily\large
        \noindent
        \@author
        \par}
  \par
  \vskip 4em}
\def\@seccntformat#1{\csname the#1\endcsname{.\ }}
\renewcommand\section{\@startsection {section}{1}{\z@}%
                                   {-3.0ex \@plus -1ex \@minus -.2ex}%
                                   {1.5ex \@plus.2ex}%
                                   {\normalfont\sffamily\large\bfseries}}
\renewcommand\subsection{\@startsection{subsection}{2}{\z@}%
                                     {-2.75ex\@plus -1ex \@minus -.2ex}%
                                     {1.5ex \@plus .2ex}%
                                   {\normalfont\sffamily\large}}
\def\fnum@figure{\normalfont\footnotesize\figurename~\thefigure}
\renewcommand\tableofcontents{%
    \section*{\contentsname
        \@mkboth{%
           \MakeUppercase\contentsname}{\MakeUppercase\contentsname}}%
    \@starttoc{toc}%
    }
\renewcommand*\l@part[2]{%
  \ifnum \c@tocdepth >-2\relax
    \addpenalty\@secpenalty
    \addvspace{2.25em \@plus\p@}%
    \begingroup
      \setlength\@tempdima{3em}%
      \parindent \z@ \rightskip \@pnumwidth
      \parfillskip -\@pnumwidth
      {\leavevmode
       \large \bfseries #1\hfil \hb@xt@\@pnumwidth{\hss #2}}\par
       \nobreak
       \if@compatibility
         \global\@nobreaktrue
         \everypar{\global\@nobreakfalse\everypar{}}%
      \fi
    \endgroup
  \fi}
\renewcommand*\l@section[2]{%
  \ifnum \c@tocdepth >\z@
    \addpenalty\@secpenalty
    \addvspace{1.0em \@plus\p@}%
    \setlength\@tempdima{1.5em}%
    \begingroup
      \parindent \z@ \rightskip \@pnumwidth
      \parfillskip -\@pnumwidth
      \leavevmode \sffamily\bfseries
      \advance\leftskip\@tempdima
      \hskip -\leftskip
      #1\nobreak\hfil \nobreak\hb@xt@\@pnumwidth{\hss #2}\par
    \endgroup
  \fi}
\renewcommand*\l@subsection{\sffamily\@dottedtocline{2}{1.5em}{2.3em}}
\renewcommand*\l@subsubsection{\@dottedtocline{3}{3.8em}{3.2em}}
\renewcommand*\l@paragraph{\@dottedtocline{4}{7.0em}{4.1em}}
\renewcommand*\l@subparagraph{\@dottedtocline{5}{10em}{5em}}
\makeatother

\makeatletter 
\@addtoreset{equation}{section}
\makeatother

\renewcommand{\theequation}{\thesection.\arabic{equation}}

\theoremstyle{plain}

\newtheorem{corollary}[equation]{Corollary}

\newtheorem{proposition}[equation]{Proposition}

\newtheorem*{theoremA}{Theorem A}
\newtheorem*{theoremB}{Theorem B}
\newtheorem*{theoremC}{Theorem C}

\theoremstyle{definition}
\newtheorem{definition}[equation]{Definition}

\newenvironment{slist}[1]%
  {\begin{list}{}%
    {%
    \settowidth{\labelwidth}{#1}%
    \setlength{\itemindent}{0pt}%
    \setlength{\labelsep}{1em}%
    \setlength{\leftmargin}{\labelwidth+\parindent+\labelsep}%
    \setlength{\itemsep}{0pt}%
    \setlength{\parsep}{.6ex}}}%
  {\end{list}}

  {\begin{list}{}%
    {%
    \settowidth{\labelwidth}{#1}%
    \setlength{\itemindent}{0pt}%
    \setlength{\labelsep}{1em}%
    \setlength{\leftmargin}{\labelwidth+\labelsep}%
    \setlength{\itemsep}{0pt}%
    \setlength{\parsep}{.6ex}}}%
  {\end{list}}

\makeatletter
\newenvironment{subequations*}{% same thing but without incrementing
  \begingroup % conservative approach
  \let\protect\@nx
  \edef\@tempa{\def\@nx\theparentequation{\theequation}}%
  \@xp\endgroup\@tempa
  \setcounter{parentequation}{\value{equation}}%
  \setcounter{equation}{0}%
  \def\theequation{\theparentequation\alph{equation}}%
  \ignorespaces
}{%
  \setcounter{equation}{\value{parentequation}}%
  \global\@ignoretrue
}

\newcommand{\eval}[2][\right]{\relax
  \ifx#1\right\relax \left.\fi#2#1\rvert}
\makeatother

\newcommand{\hide}[1]{}
\newcommand{\sdef}[1]{{\bf #1}}

\newcommand{\xra}{\xrightarrow}
\newcommand{\aut}[1]{#1}

\def\Q{\mathbb{Q}}

\def\C{\mathbb{C}}

\def\cp{\mathbb{CP}}

\newcommand{\I}{{{\mathchoice {\rm 1\mskip-4mu l} {\rm 1\mskip-4mu l}
{\rm 1\mskip-4.5mu l} {\rm 1\mskip-5mu l}}}}

\newcommand{\roots}{\Delta}

\def\la#1{{\mathfrak{#1}}}
\def\g{\la g}
\def\t{\la t}

\def\gd{\la g^*}
\def\td{\la t^*}

\def\inv{^{-1}{}}

\def\hra{\hookrightarrow}

\def\onto{\twoheadrightarrow}
\def\ddt0{{\left.{\frac{d}{dt}}\right|_{t=0}}}
\def\dds0{{\left.{\frac{d}{ds}}\right|_{s=0}}}

\def\CH{\mathcal{H}}

\def\om{\omega}

\def\al{\alpha}
\def\vp{\varphi}

\def\ep{\epsilon}

\def\lb{\lambda}
\def\Lb{\Lambda}
\def\sgm{\sigma}

\def\d/{/\mspace{-6.0mu}/}

\def\<{\left\langle}
\def\>{\right\rangle}
\def\({\left(}
\def\){\right)}
\def\iso{\cong}

\def\c.{\cdot}
\def\h..{\ldots}
\edef\oo{\o}
\def\o{\circ}

\def\ot{\otimes}
\def\ti{\times}
\def\ul{\underline}

\DeclareMathOperator{\rank}{rk}

\def\Hom{\operatorname{Hom}}
\def\End{\operatorname{End}}

\def\Tr{{\operatorname{Tr}}}
\def\tr{{\operatorname{tr}}}

\def\Lie{\operatorname{Lie}}

\def\ann{\operatorname{ann}}

\newcommand{\im}{\operatorname{im}}

\def\ideal{\vartriangleleft}

\newdir{ (}{{}*!/-1ex/@^{(}}
\newdir_{ (}{{}*!/-1ex/@_{(}}

\newcommand{\lh}{\la{h}}
\newcommand{\zg}{\mug\inv(0)}
\newcommand{\zt}{\mut\inv(0)}
\newcommand{\mug}{\mu_G}
\newcommand{\mut}{\mu_T}

\newcommand{\lm}{\bigoplus_{\al\in \roots^-} L_\al}

\newcommand{\tlp}{{\textstyle\bigoplus_{\al\in \roots^+} L_\al}}
\newcommand{\tlm}{{\textstyle\bigoplus_{\al\in \roots^-} L_\al}}
\newcommand{\inn}[2]{\int_{#1}\negthickspace\negthickspace{#2}}
\newcommand{\nts}{\negthickspace}
\DeclareMathOperator{\vvert}{vert}

\newcommand{\ccup}{{\smallsmile}}
\newcommand{\dbar}{{\overline\partial}}
\newcommand{\CE}{\mathcal{E}}
\DeclareMathOperator{\indx}{index}
\DeclareMathOperator{\ch}{ch}
\DeclareMathOperator{\Td}{Td}
\theoremstyle{plain}
\newtheorem*{theoremBp}{Theorem B${}'$}
\newtheorem*{theoremBpp}{Theorem B${}'{}'$}
\newtheorem*{theoremAH}{Theorem A${}_H$}
\newtheorem*{theoremBH}{Theorem B${}_H$}
\newtheorem*{theoremCH}{Theorem C${}_H$}

\begin{document}

\title{Symplectic quotients by a nonabelian group\\ 
and by its maximal torus}
\author{Shaun Martin%
\footnote{Institute for Advanced Study, Princeton, NJ;
smartin@ias.edu; March, 1999.}}
\maketitle

\section*{Introduction}

This paper examines the relationship between the symplectic quotient
$X\d/G$ of a Hamiltonian $G$-manifold $X$, and the associated
symplectic quotient $X\d/T$, where $T\subset G$ is a maximal torus,
in the case in which $X\d/G$ is a compact manifold or
orbifold.

The three main results are: a formula expressing the rational
cohomology ring of $X\d/G$ in terms of the rational cohomology ring of
$X\d/T$; an `integration' formula, which expresses cohomology
pairings on $X\d/G$ in terms of cohomology pairings on $X\d/T$; and an
index formula, which expresses the indices of elliptic operators on
$X\d/G$ in terms of indices on $X\d/T$.

The results of this paper are complemented by the results in a
companion paper~\cite{skm:t}, in which different
techniques are used to derive formul{\ae} for cohomology pairings on
symplectic quotients $X\d/T$, where $T$ is a torus, in terms of the
$T$-fixed points of $X$. That paper also gives some applications of
the formul{\ae} proved here.

In order to state the main results of this paper, we introduce some
notation.  The symplectic
quotient $X\d/G$ is defined to be the topological quotient
$\mug\inv(0)/G$, where $\mug:X\to\Lie(G)^*$ is a moment map for the
$G$-action on $X$. A choice of maximal torus $T\subset G$ induces a
natural projection $\Lie(G)^*\onto\Lie(T)^*$, and composing with $\mug$
gives a moment map $\mut:X\to\Lie(T)^*$ for the $T$-action, with
$X\d/T:=\mut\inv(0)/T$. In most of this paper we make some additional
simplifying assumptions: we assume that both $\mug\inv(0)$ and
$\mut\inv(0)$ are compact manifolds, on which the respective $G$- and
$T$-actions are free. It follows that $X\d/G$ and $X\d/T$ are compact
symplectic manifolds. In section~\ref{sec:gen} we show how to modify
the main results when various of these assumptions are dropped.

For every \textit{weight}
$\al$ of $T$ there is a characteristic class $e(\al)\in H^2(X\d/T)$
naturally associated%
\footnote{A weight of $T$ is a homomorphism $\al:T\to S^1$. Denoting
by $\C_{(\al)}$ the representation space on which $T$ acts via this
homomorphism, we define $L_\al\to X\d/T$ to be the associated line
bundle, that is, $L_\al:=\mut\inv(0)\ti_T\C_{(\al)}$, and $e(\al)$ to
be the Euler class of $L_\al$.}  to the principal $T$-bundle
$\mut\inv(0)\to X\d/T$ (for a precise definition see
section~\ref{sec:int}).

\begin{theoremA}[Cohomology rings] There is a natural ring isomorphism
\begin{equation*}
H^*(X\d/G;\Q) \iso \frac{H^*(X\d/T;\Q)^W}{\ann(e)}.
\end{equation*}
Here $W$ denotes the Weyl group of $G$, which acts naturally on
$X\d/T$; the class $e\in H^*(X\d/T)^W$ is defined by 
$e := \prod_{\al\in\roots} e(\al)$ for $\roots$ the
roots of $G$, and $\ann(e)\ideal H^*(X\d/T;\Q)^W$ is the ideal
consisting of all $W$-invariant elements $c\in H^*(X\d/T;\Q)^W$ such that
the product $c\ccup e$ vanishes.
\end{theoremA}

There is a natural notion of a \sdef{lift} of a cohomology class on
$X\d/G$ to a class on $X\d/T$, compatible with the above isomorphism. 
The most concrete way of expressing
this%
\footnote{A more natural way of expressing the fact that $\tilde{a}$
is a lift of $a$ brings in the $G$-equivariant cohomology
$H_G^*(X)$. There are natural maps from $H^*_G(X)$ to both
$H^*(X\d/G)$ and $H^*(X\d/T)$, and $\tilde{a}$ is a lift of $a$ if
they are both images of the same class in $H_G^*(X)$.}
involves the manifold $Y:=\mug\inv(0)/T$. There is an obvious
inclusion $i:Y\hra X\d/T$ and projection $\pi:Y\to X\d/G$, and we say
$\tilde{a}\in H^*(X\d/T)$ is a \sdef{lift} of $a\in H^*(X\d/G)$ if
$\pi^*a = i^*\tilde{a}$.

\begin{theoremB}[Integration formula]
Given a cohomology class  $a\in H^*(X\d/G)$ with lift $\tilde{a}\in
H^*(X\d/T)$, then
\begin{equation*}
\inn{X\d/G}{a} = \frac{1}{|W|}
 \inn{X\d/T}{\tilde{a}\ccup e},
\end{equation*}
where $|W|$ is the order of the Weyl group of $G$, and $e$ is the
cohomology class defined in Theorem~A.
\end{theoremB}

This formula gives cohomology pairings on $X\d/G$ because the lift of
a class is compatible with cup product: given classes $a,b\in
H^*(X\d/G)$ with lifts $\tilde{a},\tilde{b}\in H^*(X\d/T)$, then
$\tilde{a}\ccup\tilde{b}$ is a lift of $a\ccup b$.

The symplectic quotient $X\d/G$ is a compact symplectic manifold, and
it can be given an almost complex structure
$J:T(X\d/G)\to T(X\d/G)$ which is compatible with its symplectic
structure in a certain sense. Using the same prescription\footnote{%
the operator $\dbar$ is the usual Cauchy-Riemann operator from
$(0,i)$-forms to $(0,i+1)$-forms---this is well-defined on an
almost-complex manifold, although $\dbar\o\dbar$ does not necessarily
vanish---see for example Griffiths and
Harris~\cite[p.~80]{gri-har:pri-alg}; the almost complex structure
combines with the symplectic structure to define a natural metric,
which we use to define $\dbar^*$.  (The complex structure can
alternatively be viewed as defining a spin${}^c$ structure, with $D$
defined as the spin${}^c$-Dirac operator, taking even spinors to odd
spinors.)}
as on a K\"ahler manifold, we can then define an elliptic operator
\begin{equation*}
D:=\dbar+\dbar^*: C^\infty\big(\Lb^{\text{even}}T^{0,1}\big) \to
C^\infty\big(\Lb^{\text{odd}}T^{0,1}\big).
\end{equation*}
Furthermore, if
$V\to X\d/G$ is a complex vector bundle, then a choice of Hermitian
connection on $V$ lets us define an elliptic operator
$
D_V: C^\infty\big(V\ot\Lb^{\text{even}}T^{0,1}\big) \to
C^\infty\big(V\ot\Lb^{\text{odd}}T^{0,1}\big).
$
While the operator $D_V$ depends on the various choices involved, its
index does not, and we have

\begin{theoremC}[Index formula]
Suppose $V\to X\d/G$ is a complex vector bundle, and $\tilde{V}\to
X\d/T$ is a lift of $V$. Then
\begin{equation*}
\indx^{X\d/G} D_V = \indx^{X\d/T} D_{\tilde{V}\ot\Lb^{\text{even}}E}
 - \indx^{X\d/T} D_{\tilde{V}\ot\Lb^{\text{odd}}E}
\end{equation*}
Here we can take $E\to X\d/T$ to equal $\bigoplus_{\al\in\roots^+}L_\al$ for
any choice $\roots^+$ of positive roots of $G$, where $L_\al\to X\d/T$
denotes the complex line bundle naturally associated to the weight $\al$.
\end{theoremC}

This formula is simpler than one would get by applying the
Atiyah-Singer index theorem to the integration formula---its proof
runs along similar lines, but brings in a result of Borel and
Hirzebruch on the Todd genus of $G/T$.

The layout of this paper is as follows. Section~\ref{sec:top} contains
the main topological result, giving a detailed description of the
topological relationship between $X\d/G$ and
$X\d/T$. Section~\ref{sec:int} uses this result, together with some
cohomological facts, to prove theorem~B, the integration formula.  In
section~\ref{sec:coh} the integration formula is combined with
Poincar{\'e} duality to prove theorem~A, and in section~\ref{sec:ind} the
index formula is proved. Section~\ref{sec:cha} is a very short section
describing some formul{\ae} for characteristic numbers of $X\d/G$, such
as the Euler characteristic and the signature. In
section~\ref{sec:gen} various generalizations of the main results are
described, including straightforward generalizations to the case when
the two symplectic quotients are compact orbifolds, to the case in
which $X\d/T$ may have singularities or be noncompact, and finally, in
a different direction, a generalization in which $T$ is replaced by a
subgroup of full rank. Finally, in
section~\ref{sec:gra}, the results of this paper are applied
to the explicit example of the Grassmannian of $k$-dimensional planes
in $\C^n$, which arises as a symplectic quotient by the group $U(k)$.
The associated symplectic quotient by the maximal torus is the
$k$-fold product $(\cp^{n-1})^k$. One result is a
presentation of the cohomology of the Grassmannian which is
different from the usual one.

\subsection*{Relation to other results}

The results of this paper, together with the companion
paper~\cite{skm:t}, have been applied by Jeffrey and
Kirwan~\cite{jef-kir:int-the}  to prove certain formul{\ae}
for cohomology pairings on moduli spaces of stable holomorphic bundles
over a Riemann surface. These formul{\ae} were first derived by
Witten~\cite{wit:two-dim-rev} using physical arguments. Indeed the main
motivation for the results of this paper and its companion was to find
a purely topological proof of Witten's formul{\ae}.

This work was carried out in 1994, in Oxford, and at the Isaac Newton
Institute in Cambridge. After completing this work, I was made aware
of some related results of Ellingsrud and Str\oo
mme~\cite{ell-str:cho-rin}. The present paper intersects with theirs
in the case that $X$ is a complex vector space; in that case they have
a result closely related to theorem~A. (Their general setting
is the `Geometric invariant theory'
quotient~\cite{mum-fog-kir:geo-inv} of a vector space over an
arbitrary field, for which they calculate the Chow ring; their
techniques are completely different from those used here.)

\subsection*{Acknowledgements}

It is a pleasure to thank Simon Donaldson, Robert Purves, 
John Rawnsley and Dietmar
Salamon for illuminating discussions, and Lisa Jeffrey and Frances
Kirwan for their interest and encouragement while these results were
being developed.

\subsection*{Notation}

Fixed throughout this paper are the following:
\begin{slist}{${X}{p}{T}$}
\item[$X$] is a fixed smooth symplectic manifold (with symplectic form $\om$);
\item[$G\supset T$] is a connected compact Lie group and a fixed
maximal torus, both acting on $X$, preserving $\om$;
\item[$\g\supset \t$] are their Lie algebras;
\item[$\mu_G:X\to\gd$] is a moment map for the $G$ action on
$X$, which we assume throughout to be proper, and to have $0$ as a
regular value (our sign convention is $d\<\mug,\xi\>=\om(V(\xi),\c.),
\ \forall\xi\in\g$, where $V(\xi)\in\Gamma(TX)$ denotes the vector field
generated by the infinitesimal action of $\xi$; see the companion
paper~\cite{skm:t} for more details);
\item[$\mu_T:X\to\td$] is the corresponding moment map for the
restriction of the action to $T$ (given by composing $\mug$ with the
natural projection $\gd\onto\td$);
\item[$X\d/G, X\d/T$] denote the symplectic quotients
${\mug\inv(0)}/{G}$ and $\mut\inv(0)/T$ respectively.
\end{slist}

\section{The main topological result}
\label{sec:top}

The results of this paper all follow from one topological result,
which we prove in this section. This result is stated in terms of
certain complex line bundles on $X\d/T$:

\begin{definition}\label{def:l_al}
Let $\al$ be a weight of $T$, that is, a one-dimensional
representation, and let $\C_{(\al)}$ denote the corresponding
representation space. Then we define the line bundle $L_\al \to X\d/T$
to be the associated bundle
\begin{equation*}
\xymatrix{
L_\al:=\ \zt \ti_T \C_{(\al)}\qquad\qquad \ar[d] \\ X\d/T}
\end{equation*}
\end{definition}

We denote by $\roots$ the set of roots of $G$, that is, $\roots$ is
the set of nonzero weights which occur in the action of $T$ on
$\g\ot\C$; we fix a choice $\roots^+\subset\roots$ of positive
roots, and denote by $\roots^-$ the corresponding negative roots.

$X\d/G$ and $X\d/T$ are related by
a fibering and an inclusion:
\begin{equation*}
\xymatrix{ \zg/T\  \ar@{^{(}->}[r]^{i} \ar[d]^\pi
& \zt/T = X\d/T \\
**[l] X\d/G = \zg/G.
}
\end{equation*}
Note that $X\d/G$ and $X\d/T$ are symplectic manifolds, and hence
possess compactible almost complex structures, unique up to
homotopy~\cite[proposition 4.1]{mcd-sal:int-sym}.

\begin{proposition} \label{pr:top} 
\begin{enumerate}
\item
The vector bundle $\lm\to X\d/T$ has a section $s$, which is
transverse to the zero section, and such that the zeroset of $s$ is
the submanifold $\zg/T\subset X\d/T$. It follows that the derivative
of $s$ identifies the normal bundle
\begin{equation*}
\nu(\zg/T) \iso \eval{\tlm}_{\zg/T}.
\end{equation*}
\item
Letting $\vvert(\pi)\to\zg/T$ denote the vector bundle of tangents to
the fibres of $\pi$, we have
\begin{equation*}
\vvert(\pi) \iso \eval{\tlp}_{\zg/T}.
\end{equation*}
\item
There is a complex orientation\footnote{an almost complex structure on
a manifold or a vector bundle defines a complex orientation, and two
almost complex structures which are homotopic (through almost complex
structures) define the same complex orientation. For the definition of
complex orientation (which also involves stabilization)
see~\cite{sto:not-cob}.
}of $\mug\inv(0)/T$ such that the above
isomorphisms are isomorphisms of complex-oriented spaces and 
vector bundles, with respect to the complex orientations of $X\d/G$
and $X\d/T$ induced by their symplectic forms.
\end{enumerate}
\end{proposition}

A complex orientation induces a real orientation in a standard manner,
and for most of this paper we will only need that the above
isomorphisms are isomorphisms of real-oriented spaces and vector
bundles. We will need complex orientations for the results on
indices of elliptic operators and characteristic numbers, in
sections~\ref{sec:ind} and~\ref{sec:cha}.

\begin{proof} We prove the three statements in the proposition in
order.

\noindent\textit{1. The inclusion.}
We have the commuting triangle of maps
\begin{equation*}
\xymatrix{ X \ar[r]^{\mug} \ar[rd]_{\mut} 
&\gd \ar@{>>}[d]^p \\
& \td
}
\end{equation*}
where the projection $p$ is induced by the inclusion $T\hra G$.
Define $V\subset\gd$ by $V:=p\inv(0)$. Then $\zt=\mug\inv(V)$, and the
fact that $0\in\td$ is a regular value for $\mut$ is equivalent to the
assertion that the subspace $V$ is transverse to the map $\mug$. Note
that $\mut$ is a $T$-equivariant map and the coadjoint action of $T$
on $\gd$ preserves the subspace $V$. Moreover, given our choice of
positive and negative roots, we have
$V\iso\bigoplus_{\al\in\roots^-}\C_{(\al)}$.

The restriction of $\mug$ to $\zt$ defines a $T$-equivariant map
$\tilde{s}:\zt\to V$, and taking the quotient by $T$, then $\tilde{s}$
defines a section $s$ of the associated bundle $\zt\ti_T V\to X\d/T$.
Since $0\in\gd$ is a regular value of $\mug$ it follows that $0\in V$
is a regular value of $\tilde{s}$, and hence $s$ is transverse to the
zero section.

Finally, the identification of $V$ in terms of negative roots gives
$\zt\ti_T V\iso\tlm$.

\noindent\textit{2. The fibering.} 
We write $Z:=\mug\inv(0)$. Let $\pi_G$ and $\pi_T$ denote the projections
\begin{equation*}
\xymatrix{**[l] \mug\inv(0) =: Z \ar[r]^{\pi_T} \ar[rd]_{\pi_G}
  & Z/T \ar[d]^\pi \\
  & Z/G
}
\end{equation*}
($\pi$ was defined above the proposition). 

Consider the foliation of $Z$ given by the $G$-orbits: using a
$G$-invariant metric to take orthgonal complements, and the
infinitesimal action to identify tangents to the $G$-orbits, we have
the $G$-equivariant identification
\begin{equation*}
TZ \iso (Z\ti\g) \oplus \pi_G^* T(Z/G),
\end{equation*}
where $G$ acts on $\g$ by the adjoint action.
Restricting to the action of $T$, we can refine this identification
using the $T$-equivariant decomposition $\g\iso\t\oplus\la{v}$
\begin{equation*}
TZ \iso (Z\ti\t) \oplus (Z\ti\la{v}) \oplus \pi_G^* T(Z/G),
\end{equation*}
with $\la{v}\iso\bigoplus_{\al\in\roots^+}\C_{(\al)}$.
Identifying the $Z\ti\t$ factor as the tangents to the $T$-orbits,
and taking the quotient by $T$,
\begin{equation*}
T(Z/T) \iso (Z\ti_T\la{v}) \oplus \pi^* T(Z/G).
\end{equation*}
Hence, identifying $\la{v}$ in terms of positive roots gives 
$\vvert(\pi)\iso Z\ti_T\la{v} \iso \eval{\tlp}_{Z/T}$.

\noindent\textit{3. The orientation.}
We begin by summarizing the arguments in the final stage of the proof.
We will first describe the symplectic form of
$X\d/T$, restricted to $Z/T$, in terms of a decomposition of
the tangent bundle (equation~\eqref{eq:nzt} below). 
We then describe an almost complex structure $\tilde{J}_0$ on $X\d/T$
which is compatible with the symplectic form on $X\d/T$, and which 
has a simple description over $Z/T$. 
Finally, we show that $\tilde{J}_0$ is homotopic, through almost
complex structures, to an almost complex structure $\tilde{J}_1$ on
$X\d/T$ with respect to which $Z/T$ is an almost complex submanifold,
and such that the almost complex structures on $Z/T$ and on its normal
bundle agree with the complex orientations described in the statement
of the proposition.

\noindent\textit{3(i) The identification of the symplectic form.}
On $\g\oplus\gd$ we define the $G$-invariant symplectic form $\eta$
by using the duality pairing:
\begin{equation*}
\eta(\xi,\al):= \<\xi,\al\>,\quad\forall\xi\in\g,\al\in\gd,
\end{equation*}
and demanding that $\eta$ be skew-symmetric. Applying this definition
fibrewise defines an invariant symplectic form on the vector bundle
$Z\ti(\g\oplus\gd)\to Z$, which we will also denote by $\eta$.
Then there exists a $G$-equivariant isomorphism of symplectic vector bundles
over $Z=\mug\inv(0)$
\begin{equation}\label{eq:nz}
\eval{TX}_Z \iso (Z\ti\g) \oplus (Z\ti\gd) \oplus \pi_G^*T(X\d/G)
\end{equation}
where the symplectic form on the left is the restriction of
the symplectic form on $X$, and the symplectic form on the right
is the direct sum of the symplectic form on $(Z\ti\g) \oplus
(Z\ti\gd)$ given by $\eta$, and the pullback of the symplectic form on
$X\d/G$. 

This isomorphism is defined as follows. A choice of connection on the
principal $G$-bundle $Z\to X\d/G$ defines a `horizontal subbundle'
$\CH\subset TZ$, which is isomorphic to $\pi_G^*T(X\d/G)$ (the
isomorphism is induced by the derivative $d\pi_G$). Using the
inclusion $TZ\hra\eval{TX}_Z$ we can consider $\CH$ to be a subbundle
of $\eval{TX}_Z$, and we define $\CH^\om\subset\eval{TX}_Z$ to be the
symplectic complement to $\CH$, with respect to the restriction of the
symplectic form $\om$ on $X$. Standard calculations using the moment
map then imply (1) the restriction of $\om$ to $\CH$ equals the
pullback of the symplectic form on $X\d/G$, (2) the subbundle
$\CH^\om$ is a vector bundle complement to $\CH$, containing
$\vvert(\pi_G)\iso (Z\ti\g)$, and isomorphic to
$(Z\ti\g)\oplus(Z\ti\gd)$ (with the isomorphism given by choosing an
equivariant complement to $\vvert(\pi_G)$ and identifying this
complement with $Z\ti\gd$ via $d\mu$), and (3) the restriction of
$\om$ to $\CH^\om\iso (Z\ti\g)\oplus(Z\ti\gd)$ equals the symplectic
form $\eta$ defined above.

The same arguments, applied to  $T$ and $\mut\inv(0)$ in place
of $G$ and $Z=\mug\inv(0)$, give an analogous isomorphism to
that of equation~\eqref{eq:nz} above. Combining these two
isomorphisms, in a neighbourhood of $Z$, and arguing as in step 2 of
this proof gives an isomorphism of symplectic vector bundles
\begin{equation}\label{eq:nzt}
\eval{T(X\d/T)}_{Z/T} \iso (Z\ti_T\la{v}) \oplus (Z\ti_T \la{v}^*) \oplus \pi^*T(X\d/G)
\end{equation}
such that the symplectic form on the left is the restriction of the
symplectic form on $X\d/T$, and on the right is the direct sum of the
natural symplectic form on $(Z\ti_T\la{v}) \oplus (Z\ti_T \la{v}^*)$ defined
analogously to $\eta$, and the pullback of the symplectic form on
$X\d/G$.

\noindent\textit{3(ii) The almost complex structure $\tilde{J}_0$.}

Fix (1) an almost complex structure on $X\d/G$ which is
compatible with the symplectic form, 
and (2) a $T$-invariant positive-definite 
inner product on $\la{v}$. 

The inner product on $\la{v}$ gives a duality
isomorphism $\la{v}\xra{\iso}\la{v}^*$, which is $T$-equivariant, and which thus
descends to an isomorphism 
\begin{equation*}
\vp:Z\ti_T\la{v} \to Z\ti_T\la{v}^*.
\end{equation*}
We now define $J_0$. On the subbundle $\pi^*T(X\d/G)$ we define $J_0$
to equal the almost complex structure on $X\d/G$ which we fixed
above. On the subbundle $(Z\ti_T\la{v})\oplus(Z\ti_T\la{v}^*)$ we define $J_0$
to equal $\begin{pmatrix}0&-\vp\\ \vp&0\end{pmatrix}$.  One easily
checks that $J_0$ is compatible with the symplectic form on
$\eval{T(X\d/T)}_{Z/T}$, and it follows from standard results in
symplectic geometry that there exists an almost complex structure
$\tilde{J_0}$ on $X\d/T$ which is compatible with the symplectic form,
and whose restriction equals $J_0$ (see for example McDuff and
Salamon~\cite[proposition 4.1]{mcd-sal:int-sym}).

\noindent\textit{3(iii) The homotopy.}

Fix a choice of positive roots
$\roots^+\subset\roots$.
This choice of positive roots gives a $T$-invariant complex structure
on $\la{v}$, which descends to a complex structure on $Z\ti_T\la{v}$, which we
denote by $i_{\la{v}}$. Similarly, the negative roots define a complex
structure $i_{\la{v}^*}$ on $Z\ti_T\la{v}^*$, and we have
\begin{equation*}
i_{\la{v}^*} = \vp\o (-i_{\la{v}^*}) \o \vp\inv.
\end{equation*}
We now define $J_1$. On the subbundle $\pi^*T(X\d/G)$ we
define $J_1$ to equal the almost complex structure on
$X\d/G$ which we fixed above, and hence to agree with $J_0$. On the subbundle
$(Z\ti_T\la{v})\oplus(Z\ti_T\la{v}^*)$ 
we define $J_1$ to
equal $\begin{pmatrix}i_{\la{v}}&0\\ 0&i_{\la{v}^*}\end{pmatrix}$.

We now show that $J_0$ and $J_1$ are homotopic through almost
complex structures. Consider the complex linear transformations $j_0
:=\begin{pmatrix} 0&-1\\1&0\end{pmatrix}$, and $j_1:=\begin{pmatrix}
i&0\\0&-i\end{pmatrix} \in \End(\C\oplus\C)$. Since $j_0$ and $j_1$
are unitary matrices having the same
eigenvalues, there is a unitary matrix $g_1$ such that $j_1 = g_1 j_0
g_1\inv$. Let $g_t$, for $t\in[0,1]$, be a path of unitary matrices
with $g_0$ equal to the identity and $g_1$ the matrix we have just described.
Then $j_t:=g_t j_0 g_t\inv$ is a path of complex structures on the
real vector space underlying $\C\oplus\C$.

Tensoring with $\la{v}$, and using the isomorphism provided by $\vp$, we
can thus define a path of almost complex structures $J_t$ from
$J_0$ to $J_1$ (keeping the almost complex structure on the subbundle
$\pi^*T(X\d/G)$ fixed throughout). 

We can think of an almost complex structure over $X\d/T$ as a section
of a bundle with fibres $O(2n)/U(n)$, where $2n=\dim
X\d/T$~\cite[proposition 2.46]{mcd-sal:int-sym}. Thus $\tilde{J_0}$ is such a
section, and $J_t$ is a homotopy of sections restricted to the
submanifold $Z/T$. By the homotopy extension property, we can extend
$J_t$ to a homotopy $\tilde{J}_t$ of almost complex structures on
$X\d/T$. 

But $\tilde{J}_1$ has the property that $Z/T$ is an almost complex
submanifold, such that the  complex structure on $Z/T$ and on its
normal bundle agree with the complex structures on the vector bundles
in the proposition, hence completing the proof.
\end{proof}

\section{The integration formula}\label{sec:int}

In this section we prove the integration formula, theorem~B. We begin
by recalling some cohomological techniques needed in the proof.

\subsection{Integration over the fibre}

If $\pi:Y\to B$ is a fibre bundle with fibre $F$, such that $Y,B$ and
$F$ are  compact oriented manifolds, then \sdef{integration over the
fibre} is a map
\begin{equation*}
\pi_*:H^*(Y) \to H^{*-\dim F}(B)
\end{equation*}
satisfying
\begin{description}
\item[\textbf{(Multiplication)}]
${\displaystyle\qquad\qquad
\pi_*(\pi^*(b)\ccup a) = b\ccup\pi_*a,\qquad\forall a\in H^*(Y), b\in H^*(B),
}$
\item[\textbf{(Restriction)}] If $i:S\hra B$ is the inclusion a closed oriented
submanifold, then the following square of maps commutes
\begin{equation*}
\xymatrix{H^*(\pi\inv(S)) \ar[d]^{\pi'_*}  & H^{*}(Y) 
  \ar[l]^{i'{}^*} \ar[d]^{\pi_*} \\
H^{*-\dim F}(S) & H^{*-\dim F}(B), \ar[l]_{i^*}}
\end{equation*}
where $\pi'$ denotes the restriction of $\pi$ to $\pi\inv(S)$, and
$i'$ denotes the inclusion of $\pi\inv(S)$ in $Y$.
\item[\textbf{(Composition)}]
${\displaystyle\qquad\qquad
\int_B \pi_*a = \int_Y a,\qquad\forall a\in H^*(Y).
}$
\end{description}

\subsection{The Euler class of a vector bundle}

The second fact we need involves the Euler class of a vector bundle.
Suppose $V$ is a real oriented vector bundle
over a compact oriented manifold $Y$, and $s$ is a section of $V$
which is transverse to the zero section. Then the zeroset of $s$ is a
submanifold $S$ of $Y$, and denoting by $i:S\hra Y$ the inclusion, we have
\begin{equation*}
\int_S i^*a = \int_Y a\ccup e(V),\qquad\forall a\in H^*(Y)
\end{equation*}
where $e(V)\in H^{\rank V}(Y)$ is the Euler class of $V$. There is a
natural identification $\nu S\iso\eval{V}_S$ of the normal bundle to
$S$ with the restriction of $V$, and hence the orientation of $V$
induces an orientation on $\nu S$; we assume $S,\nu S$ and $Y$ are
oriented compatibly.

If instead $V$ is a complex vector bundle, with its natural real
orientation, then the Euler class of $V$ (thought of as a real vector
bundle) equals the top Chern class of $V$ (thought of as a complex
vector bundle). For proofs and more detailed explanations, see Bott
and Tu~\cite[section 6]{bot-tu:dif-for}.

\subsection{Proof of the integration formula}

Recall that we have maps 
$i:\mug\inv(0)/T\hra X\d/T$ and $\pi:\mug\inv(0)/T\to
X\d/G$ as defined in section~\ref{sec:top},
and we say that a class $\tilde{a}\in H^*(X\d/T)$ is a \sdef{lift} of
$a\in H^*(X\d/G)$ if $\pi^*(a)=i^*(\tilde{a})$.

\begin{definition}
Given a weight $\al$ of $T$, then we define $e(\al)$ to be the Euler class
\begin{equation*}
e(\al):=e(L_\al)\in H^*(X\d/T).
\end{equation*}
(These are the cohomology classes which appear in Theorems~A and~B).
\end{definition}

We are now ready to prove theorem~B (as stated in the introduction).
This integration formula is stated in terms of cohomology classes on
$X\d/G$ and $X\d/T$. However, such classes often arise from equivariant
cohomology classes on $X$ via the `Kirwan map'. After proving the
integration formula, we will set up notation for equivariant
cohomology and state a corollary in terms of equivariant classes.

\begin{proof}[Proof of Theorem~B]
Define the class $b\in H^*(X\d/T)$ by $b:=\prod_{\al\in\roots^+}
e(L_\al)$. Then $i^*b=e(\vvert(\pi))$, by proposition~\ref{pr:top}.
We can calculate the integral over the fibres $\pi_*(i^*b)\in
H^0(X\d/G)$ by restricting to a single fibre $G/T$, using the
restriction property of the pushforward.  By naturality of the Euler
class, $\eval{e(\vvert(\pi))}_{G/T} = e(\eval{\vvert(\pi)}_{G/T}) =
e(T(G/T))$, and so $\pi_*(i^*b) = \int_{G/T}{e(T(G/T))} = \chi(G/T) =
|W|$, using a standard identification of the Euler characteristic of
$G/T$~\cite{bor-hir:cha-cla}.

Then
\begin{equation*}
\begin{split}
\inn{X\d/G}{a}
 &= \frac{1}{|W|} \int_{\zg/T}{\pi^*a \ccup i^* b},
 \quad\text{applying the push-pull formula,} \\
 &= \frac{1}{|W|} \int_{X\d/T}{i_*i^*(\tilde{a} \ccup b)},
 \quad\text{composition of pushforwards; and $\tilde{a}$ is a lift of $a$,} \\
 &= \frac{1}{|W|} \int_{X\d/T}{%
 \nts\nts\nts\tilde{a} \ccup b \ccup\nts\nts
\prod_{\al\in\roots^-}\nts e(L_\al)},
 \quad\text{\parbox{3.5in}{properties of the Euler class, applied to the\\ 
       vector bundle identified in proposition~\ref{pr:top}(1),}} \\
 &= \frac{1}{|W|} \int_{X\d/T}{\tilde{a}\ccup e},
 \quad\text{by definition of $e$.}
\end{split}
\end{equation*}
\end{proof}

\section{The relationship between the cohomology rings}\label{sec:coh}

In this section we prove theorem~A, which relates the cohomology rings
of the  symplectic quotients $X\d/G$ and $X\d/T$. The proof involves
some standard machinery in equivariant cohomology, and a crucial
result concerning the `Kirwan map', which we begin by reviewing.
We also state a version of the integration formula in the language of
equivariant cohomology.

\subsection{Some key facts in equivariant cohomology}

The $G$-equivariant cohomology of the $G$-manifold $X$, which we
denote by $H_G^*(X)$, is defined to be the ordinary cohomology of the
\sdef{homotopy quotient}
\begin{equation*}
X_G := (EG \ti X)/G,
\end{equation*}
where $EG$ is a universal space for $G$: that is, $EG$ is contractible
and has a free $G$-action. For various facts in equivariant
cohomology,
see~\cite{ati-bot:mom-map,gre-hal-van:con-cur,mat-qui:sup-tho}.
We recall that, if $K\subset G$ is a subgroup, there is a natural
restriction map $r^G_K:H_G^*(X)\to H^*_K(X)$. 

A map of fundamental importance in symplectic geometry is the
\sdef{Kirwan map}, which gives a surjective ring homomorphism from the
equivariant cohomology of a symplectic manifold onto the ordinary
cohomology of its symplectic quotient.  Explicitly, we define the
Kirwan map
\begin{equation*}
\kappa_G: H^*_G(X) \to H^*(X\d/G),
\end{equation*}
by taking the restriction to $\mug\inv(0)$, and composing this with the
natural isomorphism \mbox{$H^*_G(\mug\inv(0))\xra{\iso} H^*(X\d/G)$}. This
natural isomorphism is defined in rational cohomology whenever the
$G$-action on $\mug\inv(0)$ is locally free, and we understand
$\kappa_G$ to only be defined when this is the case.  We denote the
analogous map for the maximal torus $T$ by $\kappa_T: H^*_T(X) \to
H^*(X\d/T)$.

We observe that, for any equivariant class $a\in H^*_G(X)$, the class
$\kappa_T\o r^G_T(a)$ is a lift of $\kappa_G(a)$ (see comments in the
proof of theorem~A for more elucidation on this point). We can thus
restate the integration formula, theorem~B, in a form which is more
natural in many applications:

\begin{corollary}[Integration formula in terms of equivariant
cohomology] 
For all $a\in H^*_G(X)$,
\begin{equation*}
\inn{X\d/G}{\kappa_G(a)} = \frac{1}{|W|}\inn{X\d/T}{\kappa_T\o r^G_T(a)\ccup e},
\end{equation*}
where
\begin{equation*}
e = \prod_{\al\in\roots} \kappa_T(e_T(\C_{(\al)})),
\end{equation*}
and $e_T$ denotes the $T$-equivariant Euler class.
\end{corollary}

\subsection{The proof of theorem~A}

Observe that the Weyl group $W$ of $G$ acts on $X\d/T$: since the
normalizer $N(T)$ of $T$ preserves $\mut\inv(0)$, the action of $N(T)$
on $\mut\inv(0)$ descends to an action of $W=N(T)/T$ on the quotient
$X\d/T$.

\begin{proof}[Proof of theorem~A]
Consider the fibre bundle $\zg/T\xra{\pi} X\d/G$. This has fibre
$G/T$, and the Weyl group $W$ acts on the fibres, covering the trivial
action on the base (this is the restriction of the $W$-action on
$X\d/T$). By a result of Borel~\cite[section 27]{bor:coh-esp}, the
pullback $\pi^*$ gives an isomorphism between the rational cohomology
of the base $X\d/G$ and the $W$-invariant cohomology of the total
space $\zg/T$.  This means there is a natural ring homomorphism
\begin{equation}
\vp:H^*(X\d/T)^W \to H^*(X\d/G)
\end{equation}
given by restriction to $\zg/T$ followed by the above identification,
which we will now show to be onto.

Applying the Borel-Hirzebruch result to the homotopy quotients $X_G$
and $X_T$, one can also recover the known fact that the restriction
$r^G_T$ gives an isomorphism with the 
$W$-invariants $r^G_T: H_G^*(X) \xra{\iso} H_T^*(X)^W$.
By naturality of the maps involved, we have
\begin{equation*}
\kappa_G = \vp \o \kappa_T \o r^G_T : H^*_G(X) \to H^*(X\d/G),
\end{equation*}
and since $\kappa_G$ is onto, it follows that $\vp$ is onto.

To prove the theorem,  we thus need to show that $\ker \vp = \ann(e)$. 
Let
$a\in H^*(X\d/T)^W$. Then
\begin{equation*}
\begin{split}
\vp(a)=0 &\iff \forall c\in H^*(X\d/T)^W,\quad 
  \inn{X\d/G} \vp(a)\ccup \vp(c) = 0,
  \quad\text{\parbox{2.5in}{Poincar\'e duality on $X\d/G$;\\ 
                            surjectivity of $\vp$,}} \\
 &\iff \forall c\in H^*(X\d/T)^W, \quad 
 \inn{X\d/T} a\ccup c\ccup e = 0,
 \quad\text{by the integration formula,} \\
 &\iff \forall d\in H^*(X\d/T), \quad 
 \inn{X\d/T} (a\ccup e)\ccup d = 0, \quad
 \text{\parbox{2.5in}{since we can average $d$ by $W$\\
                      (see below),}} \\
 &\iff a\ccup e = 0, \quad\text{by Poincar\'e duality on $X\d/T$.} \\
\end{split}
\end{equation*}
In the second-last step, we note that $W$ acts by symplectomorphisms
on $X\d/T$, hence preserves orientation and integrals, and since $a$
and $e$ are $W$-invariant, we can average $d$ by $W$ to obtain a
$W$-invariant class without changing the integral.
\end{proof}

\section{The index formula}\label{sec:ind} %Or should it be K-theory?

\def\kc{z}
\def\tkc{{\tilde z}}
\def\ind#1{\operatorname{ind}_{#1}}

We now prove the index formula, theorem~C, by applying the
Atiyah-Singer index theorem to the main topological result. 
For more details on the spinc${}^c$ Dirac
operators $D$ and $D_V$ (described in the introduction) see for
example~\cite[appendix D]{law-mic:spi-geo}.

In the proof we will use $K$-theory, but only in a rudimentary way,
and we recall a few facts and set up some notation.  Given a compact
space $Y$, then $K(Y)$ is a commutative ring, whose elements are
represented by formal sums (with integer coefficients) of complex
vector bundles over $Y$. Given a vector bundle $V\to Y$, we write
$[V]\in K(Y)$ for the equivalence class it represents. Addition and
multiplication in $K(Y)$ are induced by the direct sum and tensor
product of vector bundles respectively, and these operations are
extended to formal sums of vector bundles by the usual laws of a
commutative ring.

We will use the Chern character and the Todd class of complex vector
bundles. It is a standard fact that these characteristic classes only
depend on the equivalence class $[V]\in K(Y)$ of a complex vector
bundle $V\to Y$, and that these characteristic classes can be extended
to every element of $K(Y)$ by setting
\begin{equation*}
\Td([V]-[W]) = \frac{\Td(V)}{\Td(W)},\qquad \ch([V]-[W]) = \ch(V) -
\ch(W),
\end{equation*}
for all vector bundles $V,W\to Y$ (the Todd class $\Td(V)$ is a
cohomology class of mixed degree, but it has degree-$0$ part equal to
$1\in H^0(Y)$, and it follows that $\Td(V)$ has a multiplicative
inverse in the cohomology ring).  Finally, if a vector bundle $V$ is
given as a sum of line bundles $V=\bigoplus_{1\le i\le k} L_i$ then
\begin{equation*}
\Td(V) = \prod_{1\le i\le k} \frac{c_1(L_i)}{1-\exp(-c_1(L_i))},\qquad
\ch(V) = \sum_{1\le i\le k} \exp(c_1(L_i)).
\end{equation*}

We use the maps $i:\mug\inv(0)/T\hra X\d/T$ and $\pi:\mug\inv(0)/T\to
X\d/G$ as defined in section~\ref{sec:top}, and we extend the
definition of the `lift' of a cohomology class to both vector
bundles and $K$-theory in the obvious way. Thus we say a class
$\tilde{a}\in K(X\d/T)$ is a \sdef{lift} of $a\in K(X\d/G)$ if
$\pi^*a=i^*\tilde{a}$.

\begin{proof}[Proof of theorem~C]
Fix almost complex structures on $X\d/G$ and $X\d/T$, compatible with
their respctive symplectic forms. Throughout this proof we will let
$T(X\d/G)$ and $T(X\d/T)$ denote the tangent bundles, thought of as
complex vector bundles given by these almost complex structures.

Define $E:=\tlp$, as in the statement of theorem~C. Then we have
$E^*\iso\tlm$. 
The main topological result, proposition~\ref{pr:top}, implies that 
$[T(X\d/T)]-[E\oplus E^*] \in K(X\d/T)$ is a lift of $[T(X\d/G)]\in
K(X\d/G)$. Since taking characteristic classes commutes with pullback,
it follows that
\begin{equation*}
\frac{\Td T(X\d/T)}{\Td(E)\ccup\Td(E^*)}
\end{equation*}
is a lift of $\Td T(X\d/G)$.

Define the class $b\in H^*(X\d/T)$ by $b:=\Td(E)$. Then
$i^*b=\Td(\vvert(\pi))$. Now $\pi$ is the projection of a fibre
bundle, with fibres $G/T$, which arises as the global quotient of a
principal $G$ bundle by the maximal torus $T$. A result of Borel and
Hirzebruch asserts that in this case $\pi_*\Td(\vvert(\pi)) = 1$
\cite[sections 7.4 and 22.3]{bor-hir:cha-cla}.

Applying the Atiyah-Singer index theorem, and arguing as in the proof
of the integration formula,
\begin{equation*}
\begin{split}
\indx^{X\d/G} D_V &= \int_{X\d/G} \ch(V)\ccup\Td T(X\d/G) \\
&= \int_{X\d/T} \ch(\tilde{V})\ccup
\frac{\Td T(X\d/T)}{\Td(E)\ccup\Td(E^*)}\ccup b\ccup
\prod_{\al\in\roots^-} e(\al) \\
&= \int_{X\d/T} \ch(\tilde{V})\ccup
\frac{\Td T(X\d/T)}{\Td(E^*)}\ccup
\prod_{\al\in\roots^-} e(\al) \\
&= \int_{X\d/T} \ch(\tilde{V})\ccup
\Td T(X\d/T) \ccup
\prod_{\al\in\roots^-} (1-\exp(-e(\al))) \\
&= \int_{X\d/T} \ch(\tilde{V})\ccup
\Td T(X\d/T) \ccup
\prod_{\al\in\roots^+} (1-\exp(e(\al)))\\
\end{split}
\end{equation*}
But $[\Lb^{\text{even}}E]-[\Lb^{\text{odd}}E] = \sum_{i=0}^{\rank E}
(-1)^i [\Lb^iE] = \prod_{\al\in\roots^+} ([\ul{\C}]-[L_\al])$ hence,
applying the Chern character,
$\ch(\Lb^{\text{even}}E) - \ch(\Lb^{\text{odd}}E) =
\prod_{\al\in\roots^+}(1-\exp(e(\al)))$.

Combining these formul{\ae}, using additive and multiplicative properties of the
Chern character, gives
\begin{equation*}
\begin{split}
\indx^{X\d/G} D_V 
&= \int_{X\d/T} \ch(\tilde{V})\ccup
\Td T(X\d/T) \ccup
\big(\ch(\Lb^{\text{even}}E) - \ch(\Lb^{\text{odd}}E)\big) \\
&= \int_{X\d/T} \big(\ch(\tilde{V}\ot\Lb^{\text{even}}E)
-(\ch(\tilde{V}\ot\Lb^{\text{odd}}E)\big)\ccup
\Td T(X\d/T) \\
&=\indx^{X\d/T} D_{\tilde{V}\ot\Lb^{\text{even}}E}
 - \indx^{X\d/T} D_{\tilde{V}\ot\Lb^{\text{odd}}E}.
\end{split}
\end{equation*}
Note finally that the formula we have derived is stated in terms of a
choice of positive roots, but the proof does not depend on any
properties of that choice, and hence the result holds for
any choice of positive roots.
\end{proof}

\section{Characteristic numbers}\label{sec:cha}

Using the $K$-theoretic arguments from section~\ref{sec:ind}, it is a
simple matter to derive formul{\ae} which express various characteristic
numbers of $X\d/G$ in terms of characteristic numbers of $X\d/T$.
Recall that the tangent bundles $T(X\d/G)$ and $T(X\d/T)$ can be
considered as complex vector bundles in an essentially unique way, by
taking almost complex structures compatible with their symplectic forms.

In the proof of the index formula, we used the fact that
$[T(X\d/T)]-[E\oplus E^*] \in K(X\d/T)$ is a lift of $[T(X\d/G)]\in
K(X\d/G)$. It follows that
\begin{equation*}
\frac{c(T(X\d/T))}{c(E)\ccup c(E^*)}
\end{equation*}
is a lift of the total Chern class $c(T(X\d/G))$. Hence, applying the
integration formula, we get the following formula for the Euler
characteristic of $X\d/G$
\begin{equation}
\chi(X\d/G) = \frac{1}{|W|} \int_{X\d/T}
       c(T(X\d/T))\ccup \prod_{\al\in\roots}
\frac{e(\al)}{1+e(\al)}.
\end{equation}
(We are using the fact that the top Chern class equals the Euler
class.)

Similarly, taking the $L$-class,
\begin{equation}
\operatorname{signature}(X\d/G) = \frac{1}{|W|} \int_{X\d/T}
       L(T(X\d/T))\ccup \prod_{\al\in\roots}
\operatorname{tanh} e(\al).
\end{equation}

In general, for any `multiplicative characteristic class' $m$ (see
Hirzebruch~\cite[section 1]{hir:top-met}, or Milnor and
Stasheff~\cite[section 19]{mil-sta:cha-cla}), we have
\begin{equation}
\int_{X\d/G} m(T(X\d/G)) = \frac{1}{|W|} \int_{X\d/T}
       \frac{m(T(X\d/T))}{m(E)\ccup m(E^*)} \ccup 
\prod_{\al\in\roots} e(\al).
\end{equation}

\section{Generalizations}
\label{sec:gen}

In the previous sections of this paper we have assumed that both
$\mug\inv(0)$ and $\mut\inv(0)$ are compact manifolds on which the
respective $G$- and $T$-actions are free.
In this section we show how to remove some of these assumptions.
We will keep the assumption that $\mug$ is a proper map,
having $0$ as a regular value. From this it follows that $\mug\inv(0)$
is a compact manifold, on which the $G$-action is locally free, and
hence that $X\d/G$ is a compact symplectic orbifold.

\subsection{The case in which $\mut$ is proper and has $0$ as a regular
value}

If $\mut$ is proper and has $0$ as a regular value, then $X\d/T$ is a
compact symplectic orbifold. 
The arguments in
sections~\ref{sec:top}--\ref{sec:coh}, in which we proved the
integration formula and the formula relating the cohomology rings, go
through with straightforward modifications, which we now describe.

In the main topological result, proposition~\ref{pr:top}, the
line bundles $L_\al$, as well as the normal bundle and the fibering
must all be replaced by their orbifold equivalents.  (In the companion
paper to this one~\cite{skm:t}, appendix~A summarizes the main
topological and cohomological properties of orbifolds, orbifold vector
bundles, and orbifold fibre bundles, including describing how
integration over the fibre goes over to that case.)

The classes $e(\al)$ are well-defined rational cohomology classes, and
theorem~A extends to this case unchanged (rational Poincar{\'e} duality
holds for compact oriented orbifolds). Theorem~B must be modified
to take into account the existence of global finite stabilizers, and
becomes
\begin{theoremBp}[Integration formula]
If $\mug$ and $\mut$ are proper maps, both having $0$ as a regular
value, then for any class $a\in H^*(X\d/G)$ with lift $\tilde{a}$,
\begin{equation*}
\inn{X\d/G}{a} = \frac{1}{|W|}\c.\frac{o_T(\mut\inv(0))}{o_G(\mug\inv(0))}
 \inn{X\d/T}{\tilde{a}\ccup e},
\end{equation*}
where $o_G(Y)$ denotes the order of the maximal subgroup of $G$ which
fixes every point in $Y$, and $|W|$ and $e$ are as defined in
theorem~B.
\end{theoremBp}

\subsection{The case in which $\mut$ is proper, but does not have $0$
as a regular value}

The integration formula, theorem~B${}'$, can be generalized to the
case in which $0$ is not a regular value for $\mut$ in two different
ways. One way involves perturbing the value at which we take the
symplectic quotient by $T$, which we now describe. We will then describe
the other alternative, which makes use of compactly-supported
cohomology: that alternative can also handle the case in which $\mut$
is not compact.

A tubular neighbourhood of $\zg/T$ is an orbifold, since we have
assumed that $0$ is a regular value for $\mug$.  By assumption $0$ is
not a regular value for $\mut$, but by transversality there exist
regular values arbitrarily close to $0$.  Let $\ep\in\t^*$ be a
regular value, and consider the family of symplectic quotients
$X\d/T(p):=\mut\inv(p)/T$, as $p$ moves between $0$ and $\ep$.  If
$\ep$ is sufficiently close to $0$, then we can find diffeomorphisms
between a neighbourhood of $\zg/T$ and open sets in the quotients
$X\d/T(p)$. To do this, we note that the (orbifold) vector bundle $\tlm\to
X\d/T(0)$ with section $s$, defined in proposition~\ref{pr:top}, is
naturally defined over all quotients $X\d/T(p)$. For $\ep$
sufficiently small, the section $s$ will remain transverse to the
zero-section, and hence the zeros of $s$ on the symplectic quotients
$X\d/T(p)$ will be diffeomorphic to one another, along with tubular
neighbourhoods of these zerosets. Thus, for $\ep$ sufficiently small,
we have an injection $i':\mug\inv(0)/T\hra X\d/T(\ep)$, and the
main topological result, proposition~\ref{pr:top} applies with the
map $i'$ in place of the map $i$. 
A sufficient condition on $\ep$ is that there exist path joining $\ep$
and $0$, consisting entirely of regular values (except of course for
$0$). Note also that the notion of a
`lift' of a cohomology class from $X\d/G$ to $X\d/T(\ep)$ is
well-defined in this case.

\begin{theoremBpp}[Integration formula]
Suppose $\mug$ and $\mut$ are proper maps, and $0$ is a regular value for
$\mug$. Then for any regular value
$\ep\in\td$ sufficiently close to $0\in\td$, and 
any class $a\in H^*(X\d/G)$,  with lift $\tilde{a}\in H^*(X\d/T(\ep))$,
\begin{equation*}
\inn{X\d/G}{a} = \frac{1}{|W|}\c.\frac{o_T(\mut\inv(0))}{o_G(\mug\inv(0))}
 \inn{X\d/T(\ep)}{\tilde{a}\ccup e},
\end{equation*}
where $o_G(Y)$ denotes the order of the maximal subgroup of $G$ which
fixes every point in $Y$, and $|W|$ and $e$ are as defined in
theorem~B.
\end{theoremBpp}

\subsection{The case in which $\mut$ is not proper}

If $\mut$ is not proper but $0\in\td$ is a regular value, then $X\d/T$
is a noncompact orbifold, with (orbifold) line bundles $L_\al$, and
with $\mug\inv(0)/T$ as a compact suborbifold.

The section $s$ of the bundle $\tlm\to X\d/T$ has compactly-supported
zeroset, and hence the pair $(\tlm, s)$ possesses a \sdef{relative
Euler class}\footnote{Let $E\to Y$ be an oriented vector bundle over a
noncompact manifold $Y$, and let $s$ be a section whose zeroset is
compact. The bundle $E$ possesses a Thom class $\Phi$, and we define
the relative Euler class by $e(E,s):=s^*\Phi\in H^*_c(Y)$. This
cohomology class is an invariant of the homotopy class of $s$, through
homotopies for which the zeroset remains compact at all times. Given a
section in this homotopy class which is transverse to the zero section
of $E$, then $e(E,s)$ represents the Poincar{\'e} dual of the compact
submanifold given by the zeroset (Poincar{\'e} duality on a noncompact
manifold is an isomorphism between homology and compactly-supported
cohomology, see for example~\cite[propositions 6.24 and
6.25]{bot-tu:dif-for}).
These statements go over to orbifolds, with rational
cohomology.}
\begin{equation*}
e^- := e(\tlm,s) \in H^*_c(X\d/T),
\end{equation*}
lying in the cohomology with compact support of $X\d/T$. 
Setting $e^+:= e(\tlp)$ (the regular Euler class), then the product
$e:=e^+\ccup e^-$ lies in $H^*_c(X\d/T)$, hence for any class
$\tilde{a}\in H^*(X\d/T)$, the product $a\ccup e$ has compact support
and thus has a well-defined integral over
$X\d/T$.

With this interpretation of the class $e$, the integration formula of
theorem~B${}'$ holds as stated. Moreover, if 
$0$ is not a regular value of $\mut$, we can remove the
non-orbifold points and apply the above reasoning.

\subsection{Replacing $T$ by a full-rank subgroup}

Suppose $H\subset G$ is a connected closed subgroup which contains a
maximal torus $T$. Many of the results of this paper generalize in a
straightforward way to give relationships between $X\d/G$ and $X\d/H$.

Denoting by $\lh$ the Lie algebra of $H$, then moment map $\mu_H$ for
the $H$-action is defined analogously to $\mut$, by composing $\mug$
with the natural projection $\gd\onto\lh^*$. We assume for simplicity
that $\mu_H$ is proper and has $0$ as a regular value, and that the
$G$-action is free on $\mug\inv(0)$, and the $H$-action is free on
$\mu_H\inv(0)$ (there are obvious generalizations when these
conditions are not met, as described above).

Under the adjoint action of $H$, the Lie algebra $\g$ decomposes into
subrepresentations
\begin{equation*}
\g\iso \lh\oplus\la{e}.
\end{equation*}
This decomposition is compatible with the $T$-action, and
$\t\subset\lh$, hence a choice of positive roots for $T$ gives a
complex structure to the $H$-representation $\la{e}$.  We denote by
$\CE\to X\d/H$ the associated vector bundle $\mu_H\inv(0)\ti_H
\la{e} \to X\d/H$.

With these definitions, the main topological result generalizes in the
obvious way, with the bundle $\tlp$ replaced by $\CE$, and
$\tlm$ replaced by the dual bundle $\CE^*$. Thus, with the obvious
definition of a lift of a cohomology class or vector bundle, we have

\begin{theoremBH}[Integration formula]
Given a cohomology class  $a\in H^*(X\d/G)$ with lift $\tilde{a}\in
H^*(X\d/H)$, then
\begin{equation*}
\inn{X\d/G}{a} = \frac{|W(H)|}{|W(G)|}
 \inn{X\d/H}{\tilde{a}\ccup e(\CE\oplus \CE^*)},
\end{equation*}
where $W(H)$ and $W(G)$ are the Weyl groups of $H$ and $G$ respectively.
\end{theoremBH}

\begin{theoremCH}[Index formula]
Suppose $V\to X\d/G$ is a complex vector bundle, and $\tilde{V}\to
X\d/H$ is a lift of $V$. Then
\begin{equation*}
\indx^{X\d/G} D_V = \indx^{X\d/H} D_{\tilde{V}\ot\Lb^{\text{even}}\CE}
 - \indx^{X\d/T} D_{\tilde{V}\ot\Lb^{\text{odd}}\CE}
\end{equation*}
\end{theoremCH}

Theorem~A generalizes in a special case. Suppose $W(H)$ is a
\textit{normal} subgroup of $W(G)$. Then the quotient group
$W(G)/W(H)$ can be thought of as the relative Weyl group for $H$ in
$G$, and $X\d/H$ carries an action of this relative Weyl group.
In this case we have
\begin{theoremAH}[Cohomology rings] There is a natural ring isomorphism
\begin{equation*}
H^*(X\d/G;\Q) \iso \frac{H^*(X\d/H;\Q)^W}%
{\ann\big(e(\CE\oplus \CE^*)\big)},
\end{equation*}
where $W:=W(G)/W(H)$ is the relative Weyl group.
\end{theoremAH}

\newpage
\section{Example: The complex Grassmannian} \label{sec:gra}
\def\ol{\overline}
\def\u{\la{u}}
\def\id{\mathbb I}

This section contains a worked example: the complex Grassmannian of
$k$-planes in $\C^n$, which we denote
$G(k,n)$. We first describe the results of
applying theorems~A and~B, and we then describe the derivation of
these results in more detail.

The Grassmannian can be described as the symplectic quotient of the
set of complex matrices with $n$ rows and $k$ columns by the unitary
group
\begin{equation*}
G(k,n) \iso \Hom(\C^k,\C^n)\d/U(k),
\end{equation*}
where $g\in U(k)$ acts on a matrix $A\in\Hom(\C^k,\C^n)$ by 
$A\o g\inv$.

The associated symplectic quotient by the maximal torus  $T\subset U(k)$
turns out to be the $k$-fold product $(\cp^{n-1})^k$. Its cohomology
ring  is generated by degree-$2$ elements $\{u_1,\h..,u_k\}$, where
$u_i$ is the positive generator of the cohomology ring of the $i$-th
copy of $\cp^{n-1}$. The Weyl group of $U(k)$ is the
symmetric group on $k$ elements $S^k$, which acts by permuting the
factors in $(\cp^{n-1})^k$. The roots $\al$ of $U(k)$ can be enumerated by
pairs of positive integers $(i,j)$ with $1\le i,j\le k$ and $i\ne j$, 
and the cohomology class
corresponding to the root $(i,j)$ is the class
$e(\al)=u_j-u_i$. Hence, theorem~A states

\begin{proposition}
The cohomology ring of the Grassmannian $G(k,n)$ is given by
\begin{equation*}
H^*(G(k,n);\Q) \iso \frac{\Q[u_1,\h..,u_k]^{S_k}}%
{\<u_1^n,\h..,u_k^n\>:\prod_{i\ne j}(u_i-u_j)},
\end{equation*}
where the expression $I:e$ in the denominator denotes the
\textit{ideal quotient} of the ideal $I$ by the element $e$, that is,
the ideal consisting of all elements $b$ such that $b\c.e\in I$. 
\end{proposition}

The Grassmannian possesses a tautological vector bundle $V\to
G(k,n)$ of rank $k$, and the cohomology ring is generated by the Chern
classes of the dual bundle $V^*$. In the above description, $c_i(V^*)$
is represented by the $i$-th symmetric polynomial in the $u_j$.
The above description of the cohomology ring is quite different from the
usual description (which involves the Segre classes of $V$), and I have been
able to find neither a general algebraic proof of the equivalence of
the two descriptions, nor any reference to the above description in
the literature.

Theorem~B gives

\begin{proposition}
\begin{equation*}
\int_{G(k,n)} c_1(V^*)^{m_1}\ccup\h..\ccup c_k(V^*)^{m_k} = %
 \frac{1}{k!}
\operatorname{coeff}_{u_1^{n-1}\h.. u_k^{n-1}}%
\bigg(\sgm_1^{m_1}\h..\sgm_k^{m_k}\c.
      \prod_{i\ne j}(u_i-u_j)\bigg)
\end{equation*}
where $\sgm_i$ is the $i$-th elementary symmetric polynomial of the
$u_j$, and $\operatorname{coeff}_m(p)$ denotes the coefficient of the
monomial $m$ in the polynomial $p$.
\end{proposition}

\subsection{The construction of $G(k,n)$}

The symplectic structure on $\Hom(\C^k,\C^n)$ is the standard one for
a complex vector space with coordinates, namely, if 
$a_{ij}=x_{ij}+\sqrt{-1}y_{ij}$, for $1\le i\le n, 1\le j\le k$, then 
\begin{equation*}
\om:=\sum_{i,j} dx_{ij}\wedge dy_{ij}.
\end{equation*}
The moment map takes values in the dual of the Lie algebra of $U(k)$,
which can be identified with the space of Hermitian $k\ti k$ matrices.
It is a straightforward calculation using the definition of the moment
map to show that a moment map is given by
\begin{equation*}
\mu_{U(k)}(A)= A^*A-\I,
\end{equation*}
where $A^*:=\overline{A}^\tr$ 
(precisely, given a
skew-Hermitian matrix $\xi\in\Lie(U(k))$, then the pairing
$\<\mu_{U(k)}(A),\xi\>$ is given by
$\frac{i}{2}\Tr\big((A^*A-\I)\xi^*\big)$).

The $k$ column vectors of the matrix $A\in\Hom(\C^k,\C^n)$
define vectors $v_1,\h..,v_k\in \C^n$, and the $(i,j)$-entry of $A^*A$ is the
Hermitian inner product $(v_j,v_i)$. Hence $\mu_{U(k)}\inv(0)$
consists of the unitary $k$-frames in $\C^n$, and taking the quotient
by $U(k)$ gives the Grassmannian $G(k,n)$.

The maximal torus $T\subset U(k)$ of diagonal matrices has associated
moment map given by the diagonal entries of the matrix $A^*A-\I$, and
so a $k$-tuple $(v_1,\h..,v_k)$ lies in $\mut\inv(0)$ precisely when
each $v_i$ has length $1$. The torus $T$ equals the product $(S^1)^k$,
the factors of which rotate the vectors independently, 
hence identifying the quotient $\Hom(\C^k,\C^n)\d/T$
with $(\cp^{n-1})^k$. To calculate the classes $e(\al)$, consider the
conjugation action of the diagonal matrices on the skew-symmetric
matrices whose diagonal entries vanish (this is the complement of $\t$
in $\Lie(U(k))$. The matrix with diagonal entries $(\lb_1,\h..,\lb_k)$
acts by on the $(i,j)$-entry by $\lb_i\lb_j\inv$; the complex line
bundle constructed with this weight has Euler class $u_j-u_i$.

We now describe the tautological bundle $E\to G(k,n)$ in terms of the
symplectic quotient construction, so that we can identify the Chern
classes of its dual on the $T$-symplectic quotient $(\cp^{n-1})^k$.
A point of $G(k,n)$ is a $U(k)$-orbit of nondegenerate homomorphisms
$A:\C^k\to\C^n$, and the corresponding fibre of $E$ is
$\im(A)\subset\C^n$. Two points $A$ and $Ag\inv$ in the same
$U(k)$-orbit give different identifications of $\C^k$ with the
subspace $\im(A)$, and, taking this into account gives
\begin{equation*}
E \iso \mu\inv(0)\ti_{U(k)} \C^k_{(\text{def.})},
\end{equation*}
where $\C^k_{(\text{def.})}$ denotes the defining representation of
$U(k)$. Thus $E^*$ is constructed from the  dual of the defining
representation, and when we restrict this dual representation to the
maximal torus, it decomposes into $k$ one-dimensional representations,
which have associated line bundles on the
$(\cp^{n-1})^k$ with Euler classes $u_1,\h..,u_k$.
The identification of the Chern classes as elementary symmetric
polynomials then follows.

\bibliographystyle{alpha}

\end{document}